\documentclass[a4paper,12pt]{amsart}
\usepackage{ifthen}
\usepackage{graphicx}
\usepackage{mathrsfs}
\usepackage{color}
\nonstopmode \numberwithin{equation}{section}
\newtheorem{thm}{Theorem}[section]
\newtheorem{cor}[thm]{Corollary}
\newtheorem{lem}[thm]{Lemma}

\theoremstyle{definition}

\newenvironment{pf}[1][]{%
 \vskip 3mm
 \noindent
 \ifthenelse{\equal{#1}{}}%
  {{\slshape Proof. }}%
  {{\slshape #1.} }%
 }%
{\qed\bigskip}

\newcounter{alphabet}

\newcommand{\A}{{\mathcal A}}

\newcommand{\C}{{\mathbb C}}
\newcommand{\D}{{\mathbb D}}

\newcommand{\eK}{{\mathcal K}}
\newcommand{\eE}{{\mathcal E}}

\newcommand{\Nzero}{{{\mathbb N}_0}}

\newcommand{\R}{{\mathbb R}}

\newcommand{\es}{{\mathcal S}}

\newcommand{\co}{{\overline{\operatorname{co}}}}

\newcommand{\hol}{{\operatorname{Hol}}}

\renewcommand{\Im}{{\,\operatorname{Im}\,}}

\renewcommand{\Re}{{\,\operatorname{Re}\,}}

\newcommand{\Gauss}{{\null_2F_1}}

\newcommand{\aand}{{\quad\text{and}\quad}}
\newcommand{\cm}{{\mathcal T}}
\newcommand{\acm}{{\widetilde{\mathcal T}}}

\newcounter{minutes}
\divide\time by 60
\newcounter{hours}
\multiply\time by 60 \addtocounter{minutes}{-\time}

\begin{document}
\bibliographystyle{amsplain}
\title{On geometric properties of ratio of two hypergeometric functions}

\begin{center}
{\tiny \texttt{FILE:~\jobname .tex,
        printed: \number\year-\number\month-\number\day,
        \thehours.\ifnum\theminutes<10{0}\fi\theminutes}
}
\end{center}
\author[T. Sugawa]{Toshiyuki Sugawa}
\address{Graduate School of Information Sciences \\
Tohoku University\\
Aoba-ku, Sendai 980-8579, Japan}
\email{sugawa@math.is.tohoku.ac.jp}
\author[L.-M.~Wang]{Li-Mei Wang}
\address{School of Statistics,
University of International Business and Economics, No.~10, Huixin
Dongjie, Chaoyang District, Beijing 100029, China}
\email{wangmabel@163.com}


\begin{abstract}
R.~K\"ustner proved in his 2002 paper that
the function $w_{a,b,c}(z)=\Gauss(a+1,b;c;z)/\Gauss(a,b;c;z)$
maps the unit disk $|z|<1$ onto a domain convex in the
direction of the imaginary axis under some condition
on the real parameters $a,b,c.$
Here $\Gauss(a,b;c;z)$ stands for the Gaussian hypergeometric function.
In this paper, we study the order of convexity of $w_{a,b,c}.$
In particular, we partially solve the problem raised by the afore-mentioned
paper by K\"ustner.
\end{abstract}

\thanks{This research is supported by National Natural Science Foundation of China (No.
11901086).
}
\maketitle

\section{Introduction and main results}

For a domain $\Omega$ in the complex plane $\C$, let $\hol(\Omega)$ denote the set
of all holomorphic functions in $\Omega$.
Throughout the paper, we will denote by $\D$ the open unit disk
 $\{z\in\C: |z|<1\}$ and by
$\Lambda$ the slit domain $\C\setminus[1,+\infty)$.
Let $\A_1$ be the subset of $\hol(\D)$ consisting of functions $f$  normalized by
$f(0)=f'(0)-1=0$ and $\es$ be the subset of $\A_1$ consisting of
univalent functions on $\D.$

For a non-constant function $f\in \hol(\D)$,
the {\it order of convexity} is defined by
$$
\kappa(f):=1+\inf_{z\in \D}\Re \frac{zf''(z)}{f'(z)}\in[-\infty,1].
$$
Note that $\kappa(f)$ is affine invariant; that is $\kappa(af+b)=\kappa(f)$
for constants $a,b$ with $a\ne0.$
Similarly, we can define the {\it order of starlikeness} (with
respect to the point $f(0)$) of $f\in \hol(\D)$ by
$$
\sigma(f):=\inf_{z\in \D}\Re \frac{zf'(z)}{f(z)-f(0)}\in[-\infty,1].
$$
It is known that $f$ is \emph{convex}, i.e. $\kappa(f)\geq 0$
if and only if $f$ is univalent in $\D$ and $f(\D)$ is a convex domain;
and $f$ is \emph{starlike}, i.e. $\sigma(f)\geq 0$,
if and only if $f$ is univalent in $\D$ and $f(\D)$ is a starlike domain
with respect to the point $f(0)$.
For $\alpha <1$, a function $f\in\A_1$ is called \emph{starlike of order $\alpha$}
if $\sigma(f)\ge \alpha.$
We denote by $\es^{*}(\alpha)$ the class of those starlike functions
of order $\alpha$.
The function $s_\alpha(z)=z/(1-z)^{2(1-\alpha)}$ in $\es^{*}(\alpha)$
is extremal in many respects.
Note that $s_0$ is nothing but the Koebe function.

The \emph{convolution} (or \emph{Hadamard product}) $f*g$ of two
functions $f,g\in\A_1$ with power series
$f(z)=z+\sum_{n=2}^{\infty}a_nz^n$ and
$g(z)=z+\sum_{n=2}^{\infty}b_nz^n$ is defined by
$$
(f*g)(z):=z+\sum_{n=2}^{\infty}a_nb_nz^n.
$$
Obviously, $f*g\in \A_1$.

We now recall the notion of subordination between two holomporphic
functions $f$ and $g$ on $\D.$ We say that $f$ is subordinate to $g$
and write $f\prec g$ or $f(z)\prec g(z)$ if there exists a holomorphic
function $\omega$ on $\D$ such that $\omega(0)=0,$ $|\omega(z)|<1$
and $f(z) = g(\omega(z))$ for $z\in\D.$
Note that $f(0)=g(0)$ and $f(\D)\subset g(\D)$ if $f\prec g.$
When $g$ is univalent, the converse is also true.

Our present work is inspired by the following inclusion for $\alpha,\beta\in[1/2,1]:$
\begin{equation}\label{inclusion1}
\big\{z(f*g)'(z)/(f*g)(z): z\in \D\big\}\subset
\co \, h_{\alpha,\beta}(\D)
\quad\text{for}~ f\in\es^{*}(\alpha),\, g\in\es^{*}(\beta),
\end{equation}
where $\co$ stands for the closed convex hull and the function $h_{\alpha,\beta}$
is defined by
$$
h_{\alpha,\beta}(z)=\frac{z(s_{\alpha}*s_{\beta})'(z)}{(s_{\alpha}*s_{\beta})(z)}.
$$
The above relation (\ref{inclusion1}) is contained in the proof of \cite[Thm.~2.7, p.~56]{Rus:conv}.
K\"ustner \cite{Kustner:2002} posed the problem asking whether the following
subordination holds for $\alpha,\beta\in[1/2,1]$ or not:
\begin{equation}
\frac{z(f*g)'}{f*g}\prec h_{\alpha,\beta}=\frac{z(s_{\alpha}*s_{\beta})'}{s_{\alpha}*s_{\beta}}
\quad\text{for}\,f\in\es^{*}(\alpha),\,
g\in\es^{*}(\beta).
\end{equation}
If $h_{\alpha,\beta}$ is convex on
the unit disk, the above subornation follows from \eqref{inclusion1}.
Therefore, it is an interesting problem
to find conditions on $\alpha$ and $\beta$ so
that the superordinate function $h_{\alpha,\beta}$ is convex on $\D$.
Note that, as is already speculated in \cite[p.~608]{Kustner:2002} by numerical computations,
$h_{\alpha, \beta}$ is not convex for certain parameters $\alpha,\beta\in[1/2,1].$
Since $s_\alpha* s_\beta$ is not expressed in terms of elementary functions in general,
we need special function techniques to attack this problem.

The {\it Gaussian
hypergeometric function} $\Gauss(a,b;c;z)$ is defined for $z\in\D$ by
the power series expansion
$$
\Gauss(a,b;c;z)=\sum_{n=0}^{\infty}\frac{(a)_n(b)_n}{(c)_nn!}z^n,
$$
where $(a)_n$ is the {\it Pochhammer symbol}, i.e., $(a)_0=1$ and
$(a)_n=a(a+1)\cdots(a+n-1)$ for $n=1,2,\cdots$. Here $a,b,c$ are
complex numbers with $c\notin -\Nzero:=\{0,-1,-2,\cdots\}$.
Note that $\Gauss(a,b;c;z)=\Gauss(b,a;c;z)$ by definition.
Hypergeometric functions can be
analytically continued along any path in the complex plane that avoids the
branch points $1$ and $\infty$.
In particular, they are defined on $\Lambda$ as single-valued holomorphic functions.
For more properties of the hypergeometric functions, we refer to the handbook
\cite{AbramowitzStegun:1965}.

The extremal function $s_\alpha$ of $\es^*(\alpha)$  may be
expressed in terms of hypergeometric functions as
$$
s_{\alpha}(z)=\frac{z}{(1-z)^{2(1-\alpha)}}=z\Gauss(2-2\alpha,1;1;z).
$$
A simple computation shows that
$$
(s_\alpha*s_\beta)(z)=z \Gauss(2-2\alpha, 2-2\beta; 1; z).
$$
Thus the problem can be formulated in terms of the hypergeometric functions.
In relation with this problem, K\"ustner \cite{Kustner:2002} proved that
the function
\begin{equation}\label{eq:w}
w(z)=w_{a,b,c}(z)=\frac{\Gauss(a+1,b;c;z)}{\Gauss(a,b;c;z)}
\end{equation}
maps $\D$ univalently onto a domain convex in the direction of the imaginary axis
(see Lemma \ref{lem:Kust} below).
It might be noteworthy that $s_\alpha*s_\beta$ has another expression
when $\alpha=\beta=3/4.$
Indeed, $(s_{\frac34}*s_{\frac34})(z)=(2/\pi) z \eK(z)$ and $h_{\frac34,\frac34}(z)
=\eE(z)/[(1-z)\eK(z)].$
Here $\eK(z)$ and $\eE(z)$ are the complete elliptic integrals of the first and second
kind, respectively:
$$
\eK(z)=\int_0^1\frac{dt}{\sqrt{(1-t^2)(1-zt^2)}}
\aand
\eE(z)=\int_0^1\sqrt{\frac{1-zt^2}{1-t^2}}dt.
$$

Our primary aim in this paper is to estimate the order of convexity for the
mapping $w_{a,b,c}$ and then apply it to the function $h_{\alpha, \beta}.$
For convenience, we put
$
\kappa[a,b,c]=\kappa(w_{a,b,c}).
$
As we will see in \eqref{eq:der} and Lemma \ref{lem:sym} below,
$\kappa(h_{\alpha, \beta})=\kappa[2-2\alpha, 2-2\beta,1]$ and
$\kappa[a,b,c]=\kappa[b,a,c].$
Therefore, we may assume, for instance, $a\le b$ if convenient.
Our main results are the following two theorems.

\begin{thm}\label{Thm:necessary}
$\kappa[a,b,c]=-\infty$ 
for positive real parameters $a,b$ and $c$ with
$a+b-1<c<a+b+1/2$ and $(c-a)(c-b)>0.$
In particular, the function
$\Gauss(a+1,b;c;z)/\Gauss(a,b;c;z)$ is not convex on $\D$ in this case.
\end{thm}

In particular, we have $\kappa(h_{\frac34,\frac34})=\kappa[\frac12,\frac12,1]=-\infty.$
The assumption $(c-a)(c-b)>0$ cannot be dropped in general.
Indeed, if we choose $b=c>0$ and $0<a<1,$ then
$w_{a,b,b}(z)=(1-z)^{-a-1}/(1-z)^{-a}=1/(1-z)$ and $\kappa[a,b,b]=0.$

\begin{thm}\label{Thm:sufficient}
Let $a$, $b$ and $c$ be real numbers with $0<a\le b$ and $a+b+1/2\leq c\leq 1+a.$
Then the order of convexity
of $w(z)=\Gauss(a+1,b;c;z)/\Gauss(a,b;c;z)$ satisfies the following:
$$
\kappa[a,b,c]\geq
\frac{3a-b-c}{2}+\frac{b-a}{1+a+b-c}+
\frac{ac\,[2(c-a-b)-1]}{(b+c)(1+a+b-c)}.
$$
\end{thm}

Letting $c=1$ in Theorem \ref{Thm:sufficient},
we thus obtain a lower estimate of $\kappa[a,b,1]$ as follows.

\begin{cor}\label{cor:suff}
If $0<a\leq b$ and if $a+b\leq 1/2,$ then
$$
\kappa[a,b,1]\ge
\frac{3a-b-1}{2}+\frac{b+b^2-2a^2-3ab}{(1+b)(a+b)}.
$$
In particular, $w_{a,b,1}$ is convex when
$$
(3b-1)a^2+(2b^2-5b-1)a+b(1-b^2)\ge 0.
$$
\end{cor}

In view of Theorem \ref{Thm:necessary}, the
condition $a+b\leq 1/2$ in Corollary \ref{cor:suff}
is necessary for the convexity of the function $\Gauss(a+1,b;1;z)/\Gauss(a,b;1;z)$
in the range $a+b<2.$
The picture in Figure 1 was produced by Mathematica based on Theorems \ref{Thm:necessary}
and \ref{Thm:sufficient}. Note that the white region does not necessarily
mean $\kappa[a,b,1]\le 0.$

\begin{figure}[htbp]
      \centering
      \includegraphics[width=.5\textwidth]{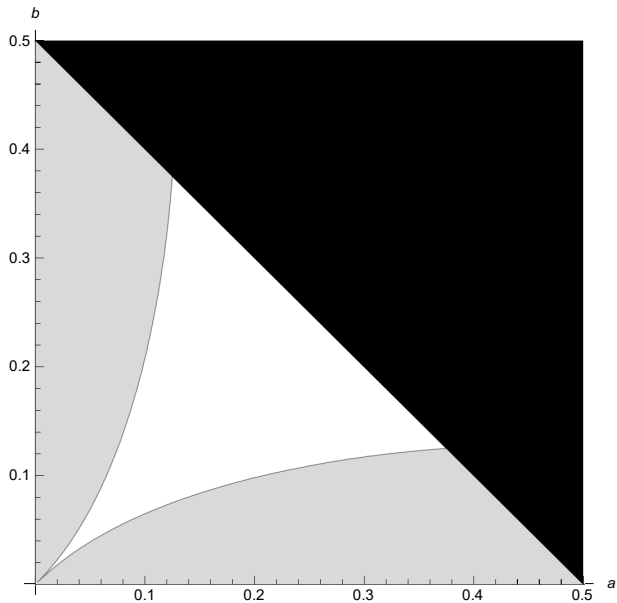}
      \caption{$\kappa[a,b,1]\ge 0$ in the gray region whereas $\kappa[a,b,1]=-\infty
$ in the black region.}
\end{figure}

Putting $a=2-2\alpha$ and $b=2-2\beta$ in Corollary \ref{cor:suff},
we arrive at the following result
solving partially the problem of K\"ustner which is mentioned above.
\begin{cor}
Let $\alpha,\beta\in[1/2,1)$ be real numbers with $\beta\le\alpha$
and $\alpha+\beta\geq 7/4$.
If
$$
2(5-6\beta)\alpha^2+(30\beta-17-8\beta^2)\alpha+4-7\beta-4\beta^2+4\beta^3\ge 0,
$$
then the subordination
$$
\frac{z(f*g)'}{f*g}\prec\frac{z(s_{\alpha}*s_{\beta})'}{s_{\alpha}*s_{\beta}},
$$
holds for $f\in\es^{*}(\alpha)$ and $g\in\es^{*}(\beta)$.
\end{cor}

\color{black}
\section{Preliminaries}
This section is devoted to some results on hypergeometric
functions which will be used in the proof of main theorems.
Formulas concerning the hypergeometric functions
without specific references below can be found
in \cite{AbramowitzStegun:1965} and \cite{OLBC:2010}.

\begin{lem}\label{lem:sym}
$\kappa[a,b,c]=\kappa[b,a,c].$
\end{lem}

\begin{proof}
By the derivative formula, we have (see \cite[(2.3)]{Kustner:2002})
\begin{equation}\label{eq:der}
\frac{z\Gauss'(a,b; c; z)}{\Gauss(a,b; c; z)}
=a\left(\frac{\Gauss(a+1,b; c; z)}{\Gauss(a,b; c; z)}-1\right).
\end{equation}
We thus see that the order of convexity $\kappa[a,b,c]$
of $w_{a,b,c}(z)=\Gauss(a+1,b; c; z)/\Gauss(a,b; c; z)$ is the same as that of
$z\Gauss'(a,b; c; z)/\Gauss(a,b; c; z),$ which is symmetric in $a$ and $b.$
Hence the required relation follows.
\end{proof}

A sequence $\{a_n\}_{n=0}^\infty$ of real numbers is called \emph{totally monotone}
or \emph{completely monotone}
if $\Delta^k a_n\ge 0$ for all integers $n,k\ge0.$
Here, $\Delta^ka_n$ is defined recursively by $\Delta^0 a_n=a_n, n\ge 0,$ and
$\Delta^{k} a_{n}=\Delta^{k-1} a_{n}-\Delta^{k-1} a_{n+1}$ for $n\geq 0$ and $k\geq 1.$

The following lemma is useful for our aim
(see Theorems 69.2 and 71.1 in \cite{Wall:anal}).

\begin{lem}\label{lem:Wall}
Let $h\in\hol(\D)$ with $h(0)\in(0,+\infty).$
Then the following three conditions are mutually equivalent:
\begin{enumerate}
\item[(i)]
$\displaystyle
h(z)=\sum_{n=0}^\infty a_nz^n,~ z\in\D,$
for a totally monotone sequence $\{a_n\}.$
\item[(ii)]
$\displaystyle
h(z)=\int_{0}^{1}\frac{d\mu(t)}{1-tz},~ z\in\D,$
for a positive Borel measure $\mu$ on $[0,1].$ 
\item[(iii)]
There exists a sequence $\{g_n\}_{n=0}^\infty$ with
$0\le g_n\le 1$ for $n=0, 1,2,\cdots$ such that
$$
h(z)=\frac{h(0)}{1-
\dfrac{(1-g_0)g_1z}{1-\dfrac{(1-g_1)g_2z}{1-\dfrac{(1-g_2)g_3z}{1-\ddots}}}},\quad
z\in\D.
$$
\end{enumerate}
\end{lem}

We denote by $\cm$ the class of functions $h\in\hol(\D)$ with $h(0)>0$
satisfying one (and hence all) of the above three conditions.
In what follows, we write
$$
h(1^-):=\lim_{x\to 1^-}h(x)=\sum_{n=0}^\infty a_n\in[0,+\infty].
$$
It is easily seen that $s_1h_1+s_2h_2\in\cm$ for $s_1, s_2\in (0,+\infty)$
and $h_1, h_2\in\cm.$
Note that a function $h\in\cm$ can be analytically continued on the domain
$\Lambda=\C\setminus[1,+\infty)$
by using the integral representation in condition (ii).
Therefore, we can regard $\cm$ as a subset of $\hol(\Lambda).$
Let $\acm$ be the class of functions $f\in\hol(\Lambda)$ of
the form $f(z)=zh(z)$ for some $h\in\cm.$
Since the sequence $a_1, a_2,\dots$ is totally monotone
for a totally monotone sequence $a_0, a_1, \dots,$
the function $h-h(0)$ belongs to $\acm$ for $h\in\cm,$ provided that $h'(0)\ne 0.$
Geometric properties of functions $f\in\acm$
are investigated by many authors (see \cite{LSW}, \cite{RusSalSu} and
\cite{Wirths75} for instance).
Among others, Wirths \cite{Wirths75} showed that a function $f$ in $\acm$ maps both
the half-plane $\Re z<1$ and the unit disk $|z|<1$
univalently onto domains convex in the direction of the imaginary axis
(see also \cite[Lem.~3.1]{Kustner:2002}).
K\"ustner \cite{Kustner:2002} studied hypergeometric functions
in connection with the class $\acm.$
Especially, the following result is important in the present work.

\begin{lem}[K\"ustner $\text{\cite[Thm.~1.5]{Kustner:2002}}$]\label{lem:Kust}
Let $a,b,c\in\R$ with $-1<a\le c$ and $0<b\le c.$
Then
$$
w_{a,b,c}(z)=\dfrac{1}{1-\dfrac{(1-g_0)g_1z}{1-\dfrac{(1-g_1)g_2z}{1-\ddots}}},
\quad z\in\D,
$$
where $g_0=0$ and
\begin{eqnarray}\label{coef}
g_n=\begin{cases}
\vspace{3mm}
\dfrac{a+k}{c+2k-1} &\quad\text{for}\,\,n=2k\ge 2,\\
\dfrac{b+k-1}{c+2k-2} &\quad\text{for}\,\,n=2k-1.
\end{cases}
\end{eqnarray}
In particular, $w_{a,b,c}\in \cm$ and the function $w_{a,b,c}$
maps the unit disk $|z|<1$ univalently onto a domain convex in the direction
of the imaginary axis.
\end{lem}

In particular, we note $g_1=b/c$ and $g_2=(a+1)/(c+1).$
We need also a few more results on continued fractions.

\begin{lem}[$\text{\cite[Thm. 11.1]{Wall:anal}}$]\label{range}
Let $g_1,g_2,g_3,\cdots$ be a sequence of real numbers which satisfies
either $0\le g_n<1~ (n=1,2,3,\dots)$ or 
$0<g_n\le 1~ (n=1,2,3,\cdots).$
Then the continued fraction
$$
v=\frac{g_1}{1-\dfrac{(1-g_1)g_2z}{1-\dfrac{(1-g_2)g_3z}{1-\ddots}}}
$$
converges uniformly on $|z|\le 1$.
Moreover, its values lie in the closed disk
$$
\left|v-\frac{1}{2-g_1}\right|\le \frac{1-g_1}{2-g_1}.
$$
\end{lem}

This estimate is simple and useful but it does not work for
a function of the form in condition (iii) of Lemma \ref{lem:Wall}
with $g_0=0.$
We offer a variant of the above estimate as follows.

\begin{lem}\label{lem:w1}
Let $h\in\cm$ with $h(1^-)<+\infty.$
Then for $z\in\D,$
$$
\left|h(z)-\frac{h(1^-)+h(-1)}{2}\right|
\le\frac{h(1^-)-h(-1)}{2}.
$$
\end{lem}

\begin{pf}
Let $\mu$ be a positive Borel measure on $[0,1]$
as in (ii) of Lemma \ref{lem:Wall}.
Noting the inequality $|z-t|\le|1-tz|$ for $|t|\le 1$ and $z\in\D,$
we now estimate
\begin{align*}
&\left|h(z)-\frac{h(1^-)+h(-1)}{2}\right|
=\left|\int_{0}^{1}\frac{d\mu(t)}{1-tz}
-\int_{0}^{1}\frac{d\mu(t)}{1-t^2}\right| \\
=\, &\left|\int_{0}^{1}\frac{t(z-t)d\mu(t)}{(1-tz)(1-t^2)}\right|
\leq\int_{0}^{1}\frac{t d\mu(t)}{1-t^2}
=\frac{h(1^-)-h(-1)}{2},\quad z\in\D.
\end{align*}
\end{pf}

The above inequality means that the value $h(z)=a_0+a_1z+a_2z^2+\cdots$ lies in the
closed disk centered at $(h(1^-)+h(-1))/2=\sum_0^\infty a_{2n}$
with radius $(h(1^-)-h(-1))/2=\sum_0^\infty a_{2n+1}.$
This refines the obvious estimate $h(-1)\le \Re h(z)\le h(1^-)$
for $z\in\D.$
When $h$ is not a constant, we can exclude the equality case,
by the maximum modulus principle.

We note here transformations of continued fractions.

\begin{lem}$($\cite[Thm.~69.2 and (75.6)]{Wall:anal}$)$
\label{transform}
Let $g_1,g_2,g_3,\cdots$ be a sequence satisfying
$0\le g_n\le 1$ for $n=1,2,3,\cdots$
and
$$
f(z)=\frac{g_1}{1-\dfrac{(1-g_1)g_2z}{1-\dfrac{(1-g_2)g_3z}{1-\ddots}}}.
$$
Then
$$
\frac{1-f(z)}{1-zf(z)}
=\frac{1-g_1}{1-\dfrac{g_1(1-g_2)z}{1-\dfrac{g_2(1-g_3)z}{1-\ddots}}},
\quad |z|\le 1,
$$
and therefore this function belongs to the class $\cm.$
\end{lem}

Applying the above lemmas, we have the following result.

\begin{lem}\label{integral2}
For $-1\le a\le c$ and $0<b\le c$, one has the expressions
$$
w(z)=\frac{\Gauss(a+1,b;c;z)}{\Gauss(a,b;c;z)}
=\frac{1}{1-zT(z)}
=\frac{1-zU(z)}{1-z}.
$$
Here, the functions $T(z)$ and $U(z)$ belong to the class $\cm$ and
they are described by
$$
T(z)=\frac{1-1/w(z)}{z}
=\frac{g_1}{1-\dfrac{(1-g_1)g_2z}{1-\dfrac{(1-g_2)g_3z}{1-\ddots}}}
$$
and
$$
U(z)=\frac{1-(1-z)w(z)}{z}
=\frac{1-T(z)}{1-zT(z)}
=\frac{1-g_1}{1-\dfrac{g_1(1-g_2)z}{1-\dfrac{g_2(1-g_3)z}{1-\ddots}}}
$$
for $|z|\le 1$, where $\{g_n\}$ is given in \eqref{coef}.
\end{lem}

Recall $g_1=b/c$ in \eqref{coef}.
Applying Lemma \ref{range} to the functions $T$ and $U,$ we obtain the following.

\begin{lem}
Let $-1\le a\le c$ and $0<b\le c.$
Then, for $|z|\le 1,$
$$
\left|T(z)-\frac{c}{2c-b}\right|\le\frac{c-b}{2c-b}
\aand
\left|U(z)-\frac{c}{b+c}\right|\le\frac{b}{b+c}.
$$
In particular,
\begin{equation}\label{eq:Re}
\frac{b}{2c-b}\le\Re T(z)\le 1
\aand
\frac{c-b}{c+b}\le\Re U(z)\le 1.
\end{equation}
\end{lem}

By the form of $T(z)$ in Lemma \ref{integral2},
we have $T(-1)=g_1/(1+(1-g_1)g_2/(1+\ddots))\le g_1=b/c.$
Using the relation $w(z)=1/(1-zT(z)),$ we have the following result
(\cite[Lem. 2.2]{Wang:2021}, see also \cite[Lem. 2.5]{Wang:20212}).

\begin{lem}\label{lem:z=-1}
If $-1\leq a\leq c$ and $0< b\leq c$, then
$$
\frac{c}{b+c}\leq w(-1)=\frac{\Gauss(a+1, b; c; -1)}{\Gauss(a, b; c; -1)}
\leq \frac{2c-b}{2c}.
$$
\end{lem}

We also need the following information about the asymptotic  behavior of $w(z)$
as $z\to 1.$

\begin{lem}\label{lem:G/F}
Let $a, b$ and $c$ be real numbers for which
none of $a,b,c,c-a,c-b$ belongs to $-\Nzero$.
The asymptotic behavior of the function
$w(z)=\Gauss(a+1,b;c;z)/\Gauss(a,b;c;z)$ as $z\to 1$ in
$\D$ $($unrestricted approach$)$ is described as follows.
\begin{enumerate}
\item[(i)]
If $a+b<c<a+b+1$, then
\begin{equation*}
w(z)=\frac{\lambda}{(1-z)^{1-\gamma}}+O(|1-z|^{\varepsilon})
\end{equation*}
where $\gamma=c-a-b\in(0,1)$, $\varepsilon=\min\{2\gamma-1,0\}$
and
\begin{equation}\label{eq:lambda}
\lambda=\frac{\Gamma(a+b+1-c)\Gamma(c-a)\Gamma(c-b)}
{\Gamma(a+1)\Gamma(b)\Gamma(c-a-b)}.
\end{equation}

\item[(ii)]
 If $c=a+b$, then
 \begin{equation*}\label{beha2}
w(z)=\frac{1}{-a(1-z)\log(1-z)}+O\left(\log\frac{1}{|1-z|}\right).
\end{equation*}

\item[(iii)]
If $a+b-1<c<a+b,$ then
\begin{equation*}
w(z)=\frac{\gamma'}{a(1-z)}
+\frac{\eta}{(1-z)^{1-\gamma'}}
+O(|1-z|^{\varepsilon'}),
\end{equation*}
where $\gamma'=a+b-c\in(0,1),$
$\varepsilon'=\min\{2\gamma'-1,0\}$
and
\begin{equation}\label{eq:eta}
\eta=\frac{\Gamma(c+1-a-b)\Gamma(a)\Gamma(b)}
{a\,\Gamma(a+b-c)\Gamma(c-a)\Gamma(c-b)}.
\end{equation}

\end{enumerate}
\end{lem}
\begin{pf}
The first two assertions can be found in
\cite[Lem. 2.3]{Wang:2021}, see also \cite[Lem. 2.4]{Wang:20212}.
We need only to
prove the last one.
Suppose that $a+b-1<c<a+b$ so that $\gamma'=-\gamma=a+b-c\in(0,1).$
Applying the general formula on $\C\setminus((-\infty,0]\cup[1,+\infty))$
\begin{align}
\Gauss(a,b;c;z)
\label{tran}
&=\frac{\Gamma(c)\Gamma(c-a-b)}{\Gamma(c-a)\Gamma(c-b)}
\Gauss(a,b;a+b-c+1;1-z) \\
&+(1-z)^{c-a-b}
\frac{\Gamma(c)\Gamma(a+b-c)}{\Gamma(a)\Gamma(b)}
\Gauss(c-a,c-b;c-a-b+1;1-z)
\notag
\end{align}
to the functions $\Gauss(a+1,b;c;z)$
and $\Gauss(a,b;c;z)$,
we obtain
\begin{align*}
\Gauss(a,b;c;z)=(1-z)^{\gamma}
\frac{\Gamma(c)\Gamma(a+b-c)}{\Gamma(a)\Gamma(b)}
\left(1+A_0(1-z)^{\gamma'}+O(1-z) \right)
\end{align*}
and
\begin{align*}
\Gauss(a+1,b;c;z)
&=(1-z)^{\gamma-1}
\frac{\Gamma(c)\Gamma(a+b-c+1)}{\Gamma(a+1)\Gamma(b)}
\left(1+O(1-z) \right) \\
&=(1-z)^{\gamma-1}\frac{a+b-c}{a}\cdot
\frac{\Gamma(c)\Gamma(a+b-c)}{\Gamma(a)\Gamma(b)}
\left(1+O(1-z) \right) ,
\end{align*}
where
$$
A_0=\frac{\Gamma(c-a-b)\Gamma(a)\Gamma(b)}%
{\Gamma(a+b-c)\Gamma(c-a)\Gamma(c-b)}
=-\frac{a\eta}{a+b-c}.
$$
Noting the functional relation $\Gamma(x+1)=x\Gamma(x),$
we obtain the required asymptotic.
\end{pf}

When $c<a+b$ and $c-a-b\ne -1,-2,\ldots$, similarly as above, we obtain
\begin{equation*}
\frac{\Gauss(a+1,b;c;z)}{\Gauss(a,b;c;z)}
=\frac{a+b-c}{a(1-z)}+O\left(|1-z|^{\gamma_1-1}\right),
\end{equation*}
where $\gamma_1=\min\{a+b-c,1\}$.
The corresponding results in
\cite[Lem. 2.3]{Wang:2021} and \cite[Lem. 2.4]{Wang:20212}
should be modified to this form.

\section{Proof of Theorem \ref{Thm:necessary}}
In this section, we will write $F(z)=\Gauss(a,b;c;z)$, $G(z)=\Gauss(a+1,b;c;z)$,
$H(z)=\Gauss(a+2,b;c;z)$,
$w(z)=G(z)/F(z)$ and $W(z)=1+zw''(z)/w'(z)$ for the sake of notational brevity.
Since the order of convexity $\kappa(w)$ of $w$ is related to the
image of $\D$ under the function $W$, we will express $W$ in terms of $w$ for a later use.
It is easy to see that $W(z)$ is expressed by
$$
W(z)=1+\frac{z(G''F-GF'')}{G'F-GF'}-\frac{2zF'}{F}.
$$
First note that the formula \eqref{eq:der} means $zF'/F=a(w-1).$
Since $F$ and $G$ satisfy the hypergeometric differential equations, we have
$$
z(1-z)F''(z)+[c-(a+b+1)z]F'(z)-ab\, F(z)=0,
$$
and
$$
z(1-z)G''(z)+[c-(a+b+2)z]G'(z)-(a+1)b\, G(z)=0.
$$
Hence, we have
$$
z(1-z)(GF''-FG'')+[c-(a+b+1)z](GF'-FG')+zFG'+bFG=0.
$$
Now we are able to rearrange the form of $W$ as
\begin{align*}
W(z)
&=1-\frac{c-(a+b+1)z}{1-z}+\frac{bFG+zG'F}{(1-z)(G'F-GF')}-2a(w-1) \\
&=\frac{1+2a-c+(b-a)z}{1-z}+\frac{b+zG'/G}{(1-z)(G'/G-F'/F)}-2aw \\
&=\frac{1+2a-c+(b-a+1)z}{1-z}+\frac{b+zF'/F}{(1-z)(G'/G-F'/F)}-2aw.
\end{align*}
By the formula $zG'/G=(a+1)(H/G-1)$ similar to \eqref{eq:der} and
Gauss' contiguous relation
$$
(a+1)(1-z)H(z)=(2a+2-c+(b-a-1)z)G(z)+(c-a-1)F(z),
$$
we compute
\begin{align*}
\frac{zG'}{G}&=\frac{(2a+2-c+(b-a-1)z)G+(c-a-1)F}{(1-z)G}-(a+1) \\
&=\frac{a+1-c+bz+(c-a-1)F/G}{1-z}
=\frac{a+1-c+bz+(c-a-1)/w}{1-z}.
\end{align*}
Finally, we arrive at the expression
\begin{equation}\label{derivative}
W(z)=a-b-1+\frac{a+b-c+2}{1-z}-2aw(z)\\
+\frac{(aw(z)+b-a)z}{Q(z)},
\end{equation}
where
\begin{equation}\label{eq:Q}
Q(z)=1+a+b-c+(a-b)(1-z)-a(1-z)w(z)+(c-a-1)/w(z).
\end{equation}
Now we are ready to prove the first main result.

\begin{pf}[Proof of Theorem \ref{Thm:necessary}]
By the definition of the order of convexity,
it suffices to prove that $\Re W(z)\to-\infty$ as $z\to 1$ along a curve in $\D.$
(If we know that $W$ is zero-free on $\D,$ we can conclude that $\kappa[a,b,c]=-\infty.$)
Indeed, we take a tangential approach to $1$ along the circle $|z-1/2|=1/2:$
$$
z_{\theta}=e^{i\theta}\cos\theta=(e^{2i\theta}+1)/2,
\quad 0<\theta<\pi/2.
$$
Note that
$$
1-z_\theta=e^{i(\theta-\pi/2)}\sin\theta
\aand
\Re\frac1{1-z_\theta}
=\frac{1-\Re z_\theta}{|1-z_\theta|^2}
=\frac{1-\cos^2\theta}{\sin^2\theta}=1.
$$
We divide the proof into three cases according to the sign of $c-a-b$.

\medskip

\noindent
{\bf Case I}: $0<\gamma=c-a-b<1/2$.
By Lemma \ref{lem:G/F} (i), we have $(1-z_\theta)w(z_\theta)
=O(\theta^\gamma)$ and
$$
w(z_\theta)=\lambda e^{i(\gamma-1)(\theta-\pi/2)}(\sin\theta)^{\gamma-1}
+O(\theta^\varepsilon) \to\infty
$$
as $\theta\to0^+,$ where $\lambda>0$ is given in \eqref{eq:lambda} and
$\varepsilon=\min\{2\gamma-1,0\}>\gamma-1.$
Hence, in view of the form of $W(z)$ in \eqref{derivative}, we get
$$
\Re W(z_\theta)=a\lambda \left(\frac{1}{1+a+b-c}-2\right)
\cos[(\gamma-1)(\theta-\pi/2)](\sin\theta)^{\gamma-1}
+O(\theta^\varepsilon).
$$
Since $\gamma<1/2,$ we conclude that $\Re W(z_\theta)\to-\infty$ as $\theta\to0^+.$

\noindent
{\bf Case II}: $a+b=c$.
By Lemma \ref{lem:G/F} (ii), we see that
 $(1-z_\theta)w(z_\theta)=O(1/(-\log\theta))$ and
\begin{align*}
w(z_\theta)
&=\frac{-1}{ae^{i(\theta-\pi/2)}\sin\theta\big[\log\sin\theta+i(\theta-\pi/2)\big]}
+O(-\log\theta) \\
&=\frac{-ie^{-i\theta}\big[\log\sin\theta-i(\theta-\pi/2)\big]}%
{a\sin\theta\big[(\log\sin\theta)^2+(\theta-\pi/2)^2\big]}+O(-\log\theta)
\to\infty
\end{align*}
as $\theta\to0^+.$
Hence,
\begin{align*}
\Re w(z_\theta)
&=\frac{-\sin\theta \log\sin\theta+(\pi/2-\theta)\cos\theta}%
{a\sin\theta\big[(\log\sin\theta)^2+(\theta-\pi/2)^2\big]}+O(-\log\theta) \\
&=\frac{\pi/2+O(-\theta\log\theta)}%
{a\theta (-\log\theta)^2}\cdot\left(1+O\left(\frac1{(-\log\theta)^{2}}\right)\right)
+O(-\log\theta) \\
&=\frac{\pi/2}{a\theta (-\log\theta)^2}+O\left(\frac1{\theta(-\log\theta)^{4}}\right)
\to+\infty
\end{align*}
as $\theta\to0^+.$
Substituting them to \eqref{derivative}, we obtain
$$
\Re W(z_\theta)=-\frac{\pi}{2\theta (-\log\theta)^2}
+O\left(\frac1{\theta(-\log\theta)^{4}}\right)
$$
and therefore $\Re W(z_\theta)\to-\infty$ as $\theta\to0^+.$

\medskip

\noindent
{\bf Case III}: $0<\gamma'=a+b-c<1.$
Lemma \ref{lem:G/F} (iii) implies that
$$
w(z_\theta)
=\frac{\gamma'}{a(1-z_\theta)}
+\eta e^{i(\gamma'-1)(\theta-\pi/2)}(\sin\theta)^{\gamma'-1}
+O(\theta^{\varepsilon'})
$$
as $\theta\to0^+,$ where $\eta$ is given in \eqref{eq:eta} and
$\varepsilon'=\min\{2\gamma'-1,0\}.$
In particular, $w(z_\theta)\to\infty$ and $a(1-z_\theta)w(z_\theta)\to
\gamma'$ as $\theta\to0^+.$
We note that $c-a>b-1>-1$ and $c-b>-1$ as well, by assumption.
Thus we see that $\eta>0$ under the present conditions.
Using \eqref{derivative}, we obtain
$$
\Re W(z_\theta)=-a\eta \cos[(\gamma'-1)(\theta-\pi/2)](\sin\theta)^{\gamma'-1}
+O(\theta^{\varepsilon'})
$$
and, in particular, $\Re W(z_\theta)\to-\infty$ as $\theta\to0^+.$
\end{pf}

We remark that the radial approach to $z=1$ does not necessarily give the
required behavior.
Indeed, for example, in Case I above, $0<\gamma<1$ and for $x\in(0,1),$
$$
W(x)=\frac{2-\gamma}{1-x}+O((1-x)^{\gamma-1})
$$
and, in particular, $\Re W(x)=W(x)\to+\infty$ as $x\to1^-.$

\section{Proof of Theorem \ref{Thm:sufficient}}

In this section, we will use the same notations as in the previous section.
We first give a slightly more general but less explicit estimate of $\kappa[a,b,c]$ than
that in Theorem \ref{Thm:sufficient}.
Recall that the functions $w=w_{a,b,c}, T$ and $U$ are defined in \eqref{eq:w} and
Lemma \ref{integral2} and they all belong to the class $\cm.$
We now define a new function $w_1$ by
\begin{equation}\label{eq:w1}
w_1(z)=a\,U(z)+(1+a-c)T(z).
\end{equation}

\begin{thm}\label{thm:tech}
Let $a$, $b$ and $c$ be real numbers with $0<a\le b$ and $a+b+1/2\leq c\leq 1+a.$
Then the order of convexity of $w=w_{a,b,c}$ is estimated as
\begin{equation}\label{eq:tech}
\kappa[a,b,c]\ge
\frac{3a-b-c}{2}+\frac{b-a}{1+a+b-c}+
\frac{a\,[2(c-a-b)-1]}{1+a+b-c}w(-1).
\end{equation}
\end{thm}

\begin{pf}
The function $Q$ defined in \eqref{eq:Q} can be described in terms of $w_1$:
\begin{eqnarray*}
Q(z)&=&1+a+b-c+(a-b)(1-z)-a(1-z)w(z)+(c-a-1)/w(z)\\
&=&1+2a-c+(b-a)z-a(1-z\, U(z))+(c-a-1)(1-zT(z))\\
&=&z \left(b-a+a\,U(z)+(1+a-c)T(z)\right)
=z(b-a+w_1(z)).
\end{eqnarray*}
Then, we can rewrite \eqref{derivative} as
\begin{equation*}
W(z)
=a-b-1+\frac{a+b-c+2}{1-z}
+\frac{b-a}{b-a+w_1(z)}+a\left(\frac{1}{b-a+w_1(z)}-2\right)w(z).
\end{equation*}
By Lemma \ref{integral2}, the function $w_1=a\,U+(1+a-c)T$ belongs to $\cm.$
In particular, $0<w_1(-1)\le w_1(1^-).$
(It is not difficult to see that equality does not hold.)
Since $w(1^-)=+\infty$ by Lemma \ref{lem:G/F}, we have $T(1^-)=1$ and
$$
U(1^-)
=1-\lim_{z\to1^-}(1-z)w(z)
=\begin{cases}
1, & \text{if}~ a+b\le c, \\
(c-b)/a, & \text{if}~ a+b>c.
\end{cases}
$$
By assumption, $c\ge a+b+1/2>a+b$ so that only the first case occurs.
Therefore,
\begin{equation}\label{eq:w1(1)}
w_1(1^-)=a\, U(1^-)+(1+a-c)T(1^-)
=1+2a-c.
\end{equation}
(Note that the above limits may be taken to be the unrestrected ones as
$z\to 1$ in $\D.$)
Now Lemma \ref{lem:w1} implies that $w_1(z)$ lies in the
closed disk $\Omega$ with the diameter $[w_1(-1), (w_1(1^-)].$
Since $w_1(-1)>0\ge a-b,$
$\Omega$ is contained in the half-plane $\{\zeta\in\C\,:\, \Re \zeta> a-b\}$.
Note that the M\"obius transformation $z\mapsto 1/(z+b-a)$
maps a closed disk in $\C\setminus\{a-b\}$ onto another closed disk.
Thus, we see that the image of $\D$ under
the mapping $1/(w_1(z)+b-a)$ lies in the disk
whose diameter is the line segment with endpoints
$1/(w_1(1^-)+b-a)$ and $1/(w_1(-1)+b-a)$.
Therefore, in conjunction with the equation \eqref{eq:w1(1)}, we have
\begin{equation}\label{eq:re2}
\Re\frac{1}{b-a+w_1(z)}\geq \frac{1}{b-a+w_1(1^-)}
=\frac{1}{1+a+b-c}\geq 2
\end{equation}
for $z\in\D$ since $1+a+b-c\leq 1/2$.
Note that $\Im z\cdot \Im h(z)\ge 0$ for $z\in\Lambda$ and $h\in\cm.$
Since $w, w_1\in\cm,$ we observe that
\begin{equation}\label{eq:im}
-\Im\left(\frac{1}{b-a+w_1(z)}-2\right)\Im w(z)
=\frac{\Im w_1(z)}{|b-a+w_1(z)|^2}\cdot\Im w(z)
\ge 0
\end{equation}
for $z\in\Lambda$.
By the inequality $\Re 1/(1-z)\ge 1/2$ for $z\in\D$ together with
\eqref{eq:re2} and \eqref{eq:im}, we have
\begin{eqnarray*}
&&\Re W(z)\\
&=&a-b-1+\Re\left[\frac{a+b-c+2}{1-z}+\frac{b-a}{b-a+w_1(z)}+
a\left(\frac{1}{b-a+w_1(z)}-2\right)w(z)\right]\\
&\geq &\frac{3a-b-c}{2}+\Re\left(\frac{b-a}{b-a+w_1(z)}\right)+
a \Re\left(\frac{1}{b-a+w_1(z)}-2\right)\Re w(z)\\
&\geq &\frac{3a-b-c}{2}+\frac{b-a}{1+a+b-c}+
a\left(\frac{1}{1+a+b-c}-2\right)w(-1).
\end{eqnarray*}
The inequality \eqref{eq:tech} now follows immediately.
\end{pf}

We note that a similar estimate may be obtained when $a>b$ in the above proof.
It seems, however, that the obtained estimate for $\kappa[a,b,c]$ is not better than
that for $\kappa[b,a,c]$ applied to \eqref{eq:tech} by interchanging $a$ and $b.$

\begin{pf}[Proof of Theorem \ref{Thm:sufficient}]
Since $w(-1)\ge c/(b+c)$ by Lemma \ref{lem:z=-1},
Theorem \ref{Thm:sufficient} now
follows from Theorem \ref{thm:tech}.
\end{pf}

We remark that the estimate in Theorem \ref{Thm:sufficient} can be
improved by using the continued fraction expansions of $w=w_{a,b,c}$
as was indicated in \cite[Rem. 2.3]{Kustner:2002}.
For instance, by Lemma \ref{lem:Kust}, we obtain
\begin{align*}
w(-1)&\ge \frac1{1+\dfrac{g_1}{1+\dfrac{(1-g_1)g_2}{1+(1-g_2)g_3}}} \\
&=\frac{(c+1) \left(-2 a b+a c+b c-2 b+c^2+4 c\right)}%
{-a b^2-2 a b c-3 a b+a c^2+a c+b^2 c+2 b c^2+3 b c+c^3+5 c^2+4 c}.
\end{align*}

\color{black}
\def\cprime{$'$} \def\cprime{$'$} \def\cprime{$'$}
\providecommand{\bysame}{\leavevmode\hbox
to3em{\hrulefill}\thinspace}
\providecommand{\MR}{\relax\ifhmode\unskip\space\fi MR }
\providecommand{\MRhref}[2]{%
  \href{http://www.ams.org/mathscinet-getitem?mr=#1}{#2}
} \providecommand{\href}[2]{#2}

\end{document}